\documentclass[12pt]{article}
\pdfoutput=1

\usepackage{a4wide}
\usepackage{amsmath}
\usepackage{amssymb}
\usepackage{amsthm}
\usepackage{picins}
\usepackage{graphicx}
\usepackage{wrapfig}
\usepackage{subfigure}
\usepackage{color}
\usepackage[small,bf]{caption}
\usepackage{ifpdf}
\ifpdf
\fi

\newtheorem{theorem}{Theorem}

\newtheorem{definition}{Definition}
\newtheorem{remark}{Remark}

\theoremstyle{definition}

\newcommand{\symm}{\bigtriangleup}

\begin{document}

\begin{center}
{\Huge On the difference between solutions of discrete tomography problems II}

\bigskip

{\Large Birgit van Dalen}

\textit{Mathematisch Instituut, Universiteit Leiden, Niels Bohrweg 1, 2333 CA Leiden, The Netherlands \\ dalen@math.leidenuniv.nl}

\today

\end{center}

{\small \textbf{Abstract:} We consider the problem of reconstructing binary images from their horizontal and vertical projections. It is known that the projections do not necessarily determine the image uniquely. In a previous paper it was shown that the symmetric difference between two solutions (binary images that satisfy the projections) is at most $4 \alpha \sqrt{2N}$. Here $N$ is the sum of the projections in one direction (i.e. the size of the image) and $\alpha$ is a parameter depending on the projections. In this paper we give a lower bound: for each set of projections that has at least two solutions, we construct two solutions that have a symmetric difference of at least $2\alpha+2$. We also show that this is the best possible.}

\section{Introduction}\label{introduction}

An important problem in discrete tomography is to reconstruct a binary image on a lattice from given projections in lattice directions \cite{boek, boeknieuw}. Each point of a binary image
has a value equal to zero or one. The line sum of a line through the image is the sum of the
values of the points on this line. The projection of the image in a certain direction
consists of all the line sums of the lines through the image in this direction.

For any set of more than two directions, the problem of reconstructing a binary image from its projections in those directions is NP-complete \cite{gardner}. For exactly two directions, the horizontal and vertical ones, say, it is possible to reconstruct an image in polynomial time. Already in 1957, Ryser described an algorithm to do so \cite{ryser}. He also characterised the set of projections that correspond to a unique binary image. Suppose $F$ is uniquely determined and has row sums $r_1$, $r_2$, \ldots, $r_m$. For each $j$ with $1 \leq j \leq \max_i r_i$ we can count the number $\#\{l: r_l \geq j\}$ of row sums that are at least $j$. Then these numbers are exactly the non-zero column sums of $F$ (in some order). See also \cite[Theorem 1.7]{boek}).

Alpers et al. \cite{alpersthesis,alpersartikel} studied the possible difference between a uniquely determined image and a second image with almost the same projections as the first one. Their results were generalised in \cite{birgit} and the same ideas were used to study the difference between two solutions of the same set of projections in \cite{birgit2}. We give an overview of the main theorems here.

Consider given row sums $\mathcal{R} = (r_1, r_2, \ldots, r_m)$ and column sums $\mathcal{C} = (c_1, c_2, \ldots, c_n)$, and assume that there exists at least one binary image with exactly these line sums. Define $\mathcal{V} = (v_1, v_2, \ldots, v_n)$ as $v_j = \#\{l: r_l \geq j\}$ for $1 \leq j \leq n$. Let $F_1$ the uniquely determined binary image with row sums $\mathcal{R}$ and column sums $\mathcal{V}$. Let $N = \sum_{i=1}^m r_i = \sum_{j=1}^n c_j$. Furthermore define the integer
\[
\alpha(\mathcal{R},\mathcal{C}) = \frac{1}{2} \sum_{j=1}^{n} |c_j - v_j|.
\]
The parameter $\alpha$ indicates how close the line sums $(\mathcal{R}, \mathcal{C})$ are to line sums that uniquely determine an image. In particular, $\alpha = 0$ if and only if there is exactly one binary image with line sums $(\mathcal{R}, \mathcal{C})$. Intuitively, the larger $\alpha$, the more possibilities there are for images that satisfy the line sums.

Alpers et al. \cite{alpersthesis,alpersartikel} proved that if $F_2$ is an image with line sums $(\mathcal{R}, \mathcal{C})$ and $\alpha(\mathcal{R}, \mathcal{C}) = 1$, then the size of the symmetric difference between $F_2$ and the uniquely determined set $F_1$ is bounded:
\[
|F_1 \symm F_2| \leq \sqrt{8N + 1} - 1.
\]
In \cite{birgit} we generalised this to larger values of $\alpha(\mathcal{R}, \mathcal{C})$. For an image $F_2$ with line sums $(\mathcal{R}, \mathcal{C})$, write $\alpha = \alpha(\mathcal{R}, \mathcal{C})$ and let $p = |F_1 \cap F_2|$. Then
\[
|F_1 \symm F_2| \leq 2\alpha + 2(\alpha+p) \log(\alpha+p)
\]
and
\[
|F_1 \symm F_2| \leq \alpha \sqrt{8N+1}-\alpha.
\]
The first bound is asymptotically sharp when $\alpha$ is large compared to $p$. The second bound is better when $\alpha$ is small compared to $p$. We have a family of examples in which this bound is achieved up to a factor $\sqrt{\alpha}$.

Using these results, we can also consider the difference between two images with the same line sums. This must be bounded by twice the upper bound for the difference between one of these images and a uniquely determined image. Hence we get the following result from \cite{birgit2}.

Suppose $F_2$ and $F_3$ are two images with the same line sums $(\mathcal{R}, \mathcal{C})$. Write $\alpha = (\mathcal{R}, \mathcal{C})$. Then
\[
|F_2 \symm F_3| \leq 2\alpha \sqrt{8N+1} - 2\alpha.
\]
As with the previous bound, there are examples for which this bound is only off by a factor $\sqrt{\alpha}$.

In this paper we consider the complementary problem: find the best lower bound for the symmetric difference between two solutions that you can at least achieve given a set of projections? For each set of projections that has at least two solutions, we construct two solutions that have a symmetric difference of at least $2\alpha+2$. We also show that this bound is sharp.

\section{Definitions and notation}\label{notation}

Let $F$ be a finite subset of $\mathbb{Z}^2$ with characteristic function $\chi$. (That is, $\chi(x,y) = 1$ if $(x,y) \in F$ and $\chi(x,y) = 0$ otherwise.) For $i \in \mathbb{Z}$, we define \emph{row} $i$ as the set $\{(x,y) \in \mathbb{Z}^2: x = i\}$. We call $i$ the index of the row. For $j \in \mathbb{Z}$, we define \emph{column} $j$ as the set $\{(x,y) \in \mathbb{Z}^2: y = j\}$. We call $j$ the index of the column. Following matrix notation, we use row numbers that increase when going downwards and column numbers that increase when going to the right.

The \emph{row sum} $r_i$ is the number of elements of $F$ in row $i$, that is $r_i = \sum_{j \in \mathbb{Z}} \chi(i,j)$. The \emph{column sum} $c_j$ of $F$ is the number of elements of $F$ in column $j$, that is $c_j = \sum_{i \in \mathbb{Z}} \chi(i,j)$. We refer to both row and column sums as the \emph{line sums} of $F$. We will usually only consider finite sequences $\mathcal{R} = (r_1, r_2, \ldots, r_m)$ and $\mathcal{C} = (c_1, c_2, \ldots, c_n)$ of row and column sums that contain all the nonzero line sums. We may assume without loss of generality that $r_1 \geq r_2 \geq \ldots \geq r_m$ and $c_1 \geq c_2 \geq \ldots \geq c_n$.

Given sequences of integers $\mathcal{R} = (r_1, r_2, \ldots, r_m)$ and $\mathcal{C} = (c_1, c_2, \ldots, c_n)$, we say that $(\mathcal{R}, \mathcal{C})$ is \emph{consistent} if there exists a set $F$ with row sums $\mathcal{R}$ and column sums $\mathcal{C}$. We say that the line sums $(\mathcal{R},\mathcal{C})$ uniquely determine such a set $F$ if
the following property holds: if $F'$ is another subset of $\mathbb{Z}^2$ with line sums $(\mathcal{R},\mathcal{C})$, then $F' = F$. In this case we call $F$ \emph{uniquely determined}.

We will now define a \emph{uniquely determined neighbour} of a set $F$. This is a uniquely determined set that is in some sense the closest to $F$. See also \cite[Section 4]{birgit2}.

\begin{definition}
Suppose $F$ has row sums $r_1 \geq r_2 \geq \ldots \geq r_m$ and column sums $c_1 \geq c_2 \geq \ldots \geq c_n$. For $1 \leq j \leq n$, let $v_j = \#\{l: r_l \geq j\}$. Then the row sums $r_1$, $r_2$, \ldots, $r_m$ and column sums $v_1$, $v_2$, \ldots, $v_n$ uniquely determine a set $F_1$, which we will call the \emph{uniquely determined neighbour of $F$}.
\end{definition}

Note that if $F'$ is another set with row sums $r_1$, $r_2$, \ldots, $r_m$ and column sums $c_1$, $c_2$, \ldots, $c_n$, then $F_1$ is a uniquely determined neighbour of $F'$ if and only if it is a uniquely determined neighbour of $F$. Hence $F_1$ only depends on the row and column sums and not on the choice of the set $F$. We will therefore also speak about the \emph{uniquely determined neighbour corresponding to the line sums} $(\mathcal{R},\mathcal{C})$, without mentioning the set $F$.

Suppose line sums $\mathcal{R} = (r_1, r_2, \ldots, r_m)$ and $\mathcal{C} = (c_1, c_2, \ldots, c_n)$ are given, where $r_1 \geq r_2 \geq \ldots \geq r_m$ and $c_1 \geq c_2 \geq \ldots \geq c_n$. Let the uniquely determined neighbour corresponding to $(\mathcal{R}, \mathcal{C})$ have column sums $v_1 \geq v_2 \geq \ldots \geq v_n$. Then we define
\[
\alpha(\mathcal{R},\mathcal{C}) = \frac{1}{2} \sum_{j=1}^{n} |c_j - v_j|.
\]
Note that $\alpha(\mathcal{R},\mathcal{C})$ is always an integer, since $2\alpha(\mathcal{R},\mathcal{C})$ is congruent to
\[
\sum_{j=1}^n (c_j + v_j) = \sum_{j=1}^n c_j + \sum_{j=1}^n v_j = 2\sum_{j=1}^n c_j \equiv 0 \mod 2.
\]

Consider a set $F$ with line sums $(\mathcal{R},\mathcal{C})$ and its uniquely determined neighbour $F_1$. Let $\alpha = \alpha(\mathcal{R},\mathcal{C})$. It was proved in \cite[Lemma 4]{birgit} that the symmetric difference $F \symm F_1$ consists of $\alpha$ \emph{staircases}. In this paper we will only use staircases of length 2, which we will define below. For the general definition of a staircase, see \cite{birgit} or \cite{birgit2}.

\begin{definition}
A \emph{staircase of length 2} in $F \symm F_1$ is a pair of points $(p_1, p_2)$ in $\mathbb{Z}^2$ such that
\begin{itemize}
\item $p_1$ and $p_2$ are in the same row,
\item $p_1$ is an element of $F \backslash F_1$,
\item $p_2$ is an element of $F_1 \backslash F$.
\end{itemize}
\end{definition}

\section{Main result}

Suppose we are given row sums $\mathcal{R} = (r_1, r_2, \ldots, r_m)$ and column sums $\mathcal{C} = (c_1, c_2, \ldots, c_n)$, where $r_1 \geq r_2 \geq \ldots \geq r_m$ and $c_1 \geq c_2 \geq \ldots \geq c_n$. Assume that the line sums are consistent but do not uniquely determine a set $F$ (hence at least two different sets with these line sums exist). Let $\alpha = \alpha(\mathcal{R},\mathcal{C})$.

In \cite{birgit2} it was shown that for all $F_2$ and $F_3$ satisfying these line sums, we have
\[
|F_2 \symm F_3| \leq 4 \alpha \sqrt{2|F_2|}.
\]
One may wonder how close we can get to achieving this bound. Our theorem shows that we can construct two sets that have a symmetric difference of size at least $2\alpha +2$.

\begin{theorem}\label{stelling}
Let row sums $\mathcal{R} = (r_1, r_2, \ldots, r_m)$ and column sums $\mathcal{C} = (c_1, c_2, \ldots, c_n)$ be given, where $r_1 \geq r_2 \geq \ldots \geq r_m$ and $c_1 \geq c_2 \geq \ldots \geq c_n$. Assume that the line sums are consistent but do not uniquely determine a set $F$. Let $\alpha = \alpha(\mathcal{R},\mathcal{C})$. Then there exist sets $F_2$ and $F_3$ with these line sums such that
\[
|F_2 \symm F_3| \geq 2\alpha + 2.
\]
This bound is sharp: for each $\alpha \geq 1$ there are line sums $(\mathcal{R},\mathcal{C})$ with $\alpha = \alpha(\mathcal{R},\mathcal{C})$ such that for any $F_2$ and $F_3$ satisfying these line sums we have $|F_2 \symm F_3| \leq 2\alpha + 2$.
\end{theorem}

\section{Proof}

In this entire section, the row sums $\mathcal{R} = (r_1, r_2, \ldots, r_m)$ and column sums $\mathcal{C} = (c_1, c_2, \ldots, c_n)$ with $r_1 \geq r_2 \geq \ldots \geq r_m$ and $c_1 \geq c_2 \geq \ldots \geq c_n$ are fixed. Furthermore, $F_1$ is the uniquely determined neighbour corresponding to $(\mathcal{R},\mathcal{C})$, and $\alpha = \alpha(\mathcal{R},\mathcal{C})$. We denote the column sums of $F_1$ by $v_1 \geq v_2 \geq \ldots \geq v_n$.

The proof is constructive. We will construct $F_2$ and $F_3$ such that they have the desired property. We will do this by changing a set $F$ step by step. Only the final result of the construction will be called $F_2$ (or $F_3$); the intermediate sets will always be called $F$ or $F'$. In Section \ref{sectionexample} the construction is illustrated by an example.

Let the columns $j$ for which $v_j > c_j$ have indices $j_1 \leq j_2 \leq \ldots \leq j_\alpha$, where each such $j$ occurs $v_j - c_j$ times. Similarly, let the columns $i$ for which $v_i < c_i$ have indices $i_1 \leq i_2 \leq \ldots \leq i_\alpha$, where each such $i$ occurs $c_i - v_i$ times. Define a \emph{column pair} as a pair $(i_t, j_t)$. The consistency of the given line sums assures that $i_t > j_t$ for all $t$. For convenience, define $i_0 = j_0 = 0$ and $i_{\alpha+1} = j_{\alpha+1} = n+1$.

We will construct both $F_2$ and $F_3$ by starting from $F=F_1$ and then for each $t$ moving an element of $F$ from column $j_t$ to column $i_t$ in the same row. After we have done that for $t = 1 , 2, \ldots, \alpha$, the row sums of $F$ have not changed, while the columns of $F$ have changed from $v_1$, $v_2$, \ldots, $v_n$ to $c_1$, $c_2$, \ldots, $c_n$. The symmetric difference $F_1 \symm F$ then consists of $\alpha$ staircases of length 2. Each staircase is confined to a single row and corresponds to a column pair $(i_t, j_t)$. We will show that we have a certain freedom in choosing the staircases.

Suppose we have moved an element for each of the column pairs $(i_1, j_1)$, $(i_2, j_2)$, \ldots, $(i_{t-1}, j_{t-1})$, where $t \geq 1$. The resulting set is called $F$ and has column sums $c_1'$, $c_2'$, \ldots, $c_n'$. Now we want to move an element from column $j_t$ to column $i_t$. For this we need a row $l$ such that the point $(l, j_t) \in F$ and $(l, i_t) \not\in F$. We have $c_{j_t}' > c_{j_t} \geq c_{i_t} > c_{i_t}'$, so $c_{j_t}' \geq c_{i_t}' + 2$. Hence there must be at least two rows that contain an element of $F$ in column $j_t$ but not in column $i_t$. This proves the existence of such a row $l$, and in fact at least two choices for $l$ are possible. Now we move the element $(l,j_t)$ to $(l, i_t)$. The row sums of $F$ do not change, while the column sum of column $j_t$ decreases by one and the column sum of column $i_t$ increases by one.

We construct both $F_2$ and $F_3$ using the construction above. First we construct $F_2$, making arbitrary choices for the rows in which we move elements. Then we will construct $F_3$. For this we let the choices in the construction depend on $F_2$, in a way we will describe below.

Let $P_1$, $P_2$, \ldots, $P_r$ be the \emph{distinct} column pairs, where $P_h$ has multiplicity $k_h$: the column pair $P_1$ is equal to each of the pairs $(i_1, j_1)$, \ldots, $(i_{k_1}, j_{k_1})$, the column pair $P_2$ is equal to each of the pairs $(i_{k_1+1}, j_{k_1+1})$, \ldots, $(i_{k_1+k_2}, j_{k_1+k_2})$, and so on. We have $k_1 + k_2 + \cdots + k_r = \alpha$. For two consecutive column pairs $(i_t, j_t)$ and $(i_{t+1}, j_{t+1})$ that are not equal we have $i_{t+1} > i_t$, $j_{t+1} \geq j_t$ or $i_{t+1} \geq i_t$, $j_{t+1} > j_t$, so the second pair contains a column that did not occur in any of the previous pairs. This means that in $P_1$, \ldots, $P_r$ at least $r+1$ different columns are involved. For each $P_h$, we denote one of the columns in $P_h$ as the \emph{final} column of $P_h$ in the following way.
\begin{itemize}
\item If one of the columns in $P_h$ also occurs in $P_{h+1}$, then the other does not occur in $P_{h+1}$, \ldots, $P_r$. We call the latter the final column of the pair.
\item If both columns in $P_h$ do not occur in $P_{h+1}$, \ldots, $P_r$, and one of the columns occurs in $P_{h-1}$, then the other does not occur in $P_1$, \ldots, $P_{h-1}$. We call the former the final column of the pair.
\item If both columns in $P_h$ do not occur in $P_1$, \ldots, $P_{h-1}$ nor in $P_{h+1}$, \ldots, $P_r$, then we arbitrarily pick one of the columns in $P_h$ and call it the final column of the pair.
\end{itemize}
By definition, we have the following properties: the final column of $P_h$ does not occur in $P_{h+1}$, \ldots, $P_r$, and if the other column does not occur in $P_{h+1}$, \ldots, $P_r$ either, then the latter column only occurs in $P_h$.

Our goal is to construct $F_3$ in such a way that, for all $h$, in the final column of $P_h$ the symmetric difference between $F_2$ and $F_3$ is at least $2k_h$, while in any other column that occurs in one of the column pairs the symmetric difference between $F_2$ and $F_3$ is at least $2$. (There is at least one such a column, since there are exactly $r$ final columns, while at least $r+1$ columns are involved in the column pairs.) If we can achieve that, then we have
\[
|F_2 \symm F_3| \geq 2k_1 + 2k_2 + \ldots + 2k_r + 2 = 2\alpha + 2.
\]

To achieve this, we choose the rows in which elements are moved for all equal column pairs at once. First we choose the rows for all pairs equal to $P_1$, then for all pairs equal to $P_2$, and so on.

Let $t$ be the index of the last column pair in a sequence of $k$ equal column pairs
\[
(i_{t-k+1},j_{t-k+1})= (i_{t-k+2}, j_{t-k+2}) =  \ldots = (i_t, j_t),
\]
where $(i_{t-k},j_{t-k}) \neq (i_{t-k+1}, j_{t-k+1})$ and $(i_t,j_t) \neq (i_{t+1}, j_{t+1})$. Suppose we have moved elements already for the column pairs $(i_1, j_1)$, \ldots, $(i_{t-k}, j_{t-k})$. Call the resulting set $F$, with column sums $c_1'$, \ldots, $c_n'$. Assume that $i_t$ is the final column of $(i_t, j_t)$ (the case where $j_t$ is the final column, is analogous). So we have $i_{t} \neq i_{t+1}$. Also, we have one of the following two properties:
\begin{itemize}
\item[(A)] $j_t = j_{t+1}$,
\item[(B)] $j_t \neq j_{t+1}$, and $j_{t-k} \neq j_{t-k+1}$.
\end{itemize}

As this is the last time column $i_t$ occurs, we need to choose the rows in such a way that by moving the elements of $F$ the symmetric difference between $F$ and $F_2$ in this column becomes at least $2k$. Also, in case (B) we want the symmetric difference in column $j_t$ to be at least 2.

Since we need to move $k$ elements out of column $j_t$ into column $i_t$, we have $c_{j_t}' \geq c_{j_t} + k \geq c_{i_t}+k \geq c_{i_t}' + 2k$, so there are at least $2k$ rows $l$ such that $(l, j_t) \in F$ and $(l, i_t) \not\in F$. Let $R$ be the set of those $2k$ rows. (If there are more than $2k$ possible rows, then pick $2k$ of them.) We distinguish between two cases.

\textit{Case 1.} Suppose there are $k$ different rows $l$ in $R$ such that $(l, i_t) \not\in F_2$. Then we move elements from column $j_t$ to column $i_t$ in each of those $k$ rows. Call the resulting set $F'$. We have $(l, i_t) \in F' \backslash F_2$ for $k$ different values of $l$. The number of elements of $F'$ in column $i_t$ must be equal to the number of elements of $F_2$ in column $i_t$, so there are also $k$ different values of $l$ for which $(l, i_t) \in F_2 \backslash F'$. Hence the symmetric difference between $F'$ and $F_2$ in this column is at least $2k$.

In case (A) we are now done, as column $j_t$ will be handled in a later column pair. Suppose we are in case (B). The column $j_t$ only occurs in the column pairs $(i_{t-k+1}, j_{t-k+1})$, \ldots, $(i_t, j_t)$, which are all equal. If for a row $l$ we have $(l, i_t) \not\in F_2$, then in the construction of $F_2$ this row was not used for a staircase corresponding to the column pair $(i_t, j_t)$ (or one of the equal ones), so we must have $(l, j_t) \in F_2$. Hence after moving elements we have $k$ different values of $l$ for which $(l, j_t) \in F_2 \backslash F'$. So in column $j_t$ the symmetric difference between $F'$ and $F_2$ is at least $2k \geq 2$.

\textit{Case 2.} Suppose there are at least $k+1$ different rows $l$ in $R$ such that $(l, i_t) \in F_2$. Let $R'$ be a set of $k+1$ of those rows. Pick one of the rows in $R'$ and call it $l_0$. Let $R''$ consist of $l_0$ and the $k-1$ other rows in $R \backslash R'$ (for which it may or may not hold that $(l, i_t) \in F_2$). Move elements from column $j_t$ to column $i_t$ in each of the $k$ rows in $R''$. Call the resulting set $F'$. Then for all $k$ rows $l$ in $R \backslash R''$ we have $(l, i_t) \in F_2 \backslash F'$. Similarly to above, we find that the symmetric difference between $F'$ and $F_2$ in column $i_t$ is at least $2k$.

Again, in case (A) we are done. Suppose we are in case (B). As column $j_t$ only occurs in the column pairs $(i_{t-k+1}, j_{t-k+1})$, \ldots, $(i_t, j_t)$, which are all equal, for at most $k$ rows $l$ in $R$ we have $(l, j_t) \not\in F_2$. This means that we can choose $l_0$ above in such a way that $(l_0, j_t) \in F_2$. After moving the elements, we then have $(l_0, j_t) \in F_2 \backslash F'$. So the symmetric difference between $F'$ and $F_2$ in column $j_t$ is at least 2.

At least one of Case 1 and Case 2 above must hold, since there are $2k$ rows in $R$. Therefore we have finished the construction of $F_2$ and $F_3$ such that $F_2 \symm F_3 \geq 2\alpha + 2$.

We will now prove the second part of Theorem \ref{stelling} by giving a family of examples for which the bound of $2 \alpha + 2$ is sharp. Let $s \geq 1$ be an integer. Take $m = n = s+1$ and let all row and column sums be equal to 1. These line sums are consistent. The uniquely determined neighbour $F_1$ has column sums $v_1 = s+1$, $v_2 = v_3 = \ldots = v_{s+1} = 0$, so $\alpha = s$.

Suppose $F_2$ and $F_3$ satisfy the given row and column sums. We have $|F_2| = |F_3| = s+1$, hence
\[
|F_2 \symm F_3| \leq |F_2| + |F_3| = 2(s+1) = 2\alpha + 2.
\]

This completes the proof of Theorem \ref{stelling}. \hfill $\square$

\begin{remark}
There do not seem to be very many examples for which the bound of $2\alpha +2$ is sharp. In particular, they all seem to have $m=n = \alpha + 1$. However, even in more general cases, when $\alpha$ is much larger than $n$, the bound is not very far off. Take for example $m=n$ and let all line sums be equal to $k$, where $k \leq \frac{1}{2}n$. The uniquely determined neighbour has $k$ column sums equal to $n$ and $n-k$ column sums equal to $0$, so $\alpha = k(n-k)$. As $n-k \geq \frac{1}{2}n$, we have $\alpha \geq \frac{1}{2}kn$. Suppose $F_2$ and $F_3$ satisfy the given row and column sums, then $|F_2| = |F_3| = kn$, hence
\[
|F_2 \symm F_3| \leq |F_2| + |F_3| = 2kn \leq 4\alpha.
\]
\end{remark}

\section{Example}\label{sectionexample}


We illustrate the construction in the proof by an example. Let row sums $(5, 5, 5, 4, 4, 2, 1, 1)$ and column sums $(6,6,6,3,3,3)$ be given. The uniquely determined neighbour $F_1$ has the same row sums, but column sums $(8,6,5,5,3,0)$ (see Figure \ref{plaatjeF1}). From this we derive that $\alpha = 4$ and that the four column pairs are $(1,3)$, $(1,6)$, $(4,6)$ and $(4,6)$.

\setlength{\intextsep}{-2pt}
\begin{wrapfigure}{R}{0pt}
\subfigure[The set $F_1$ with its row and column sums.]{\label{plaatjeF1}\includegraphics{constructie.1}}
\qquad
\subfigure[The set $F_2$ with its row and column sums.]{\label{plaatjeF2}\includegraphics{constructie.4}}
\caption{}
\label{plaatjeF2beide}
\end{wrapfigure}
\setlength{\intextsep}{3pt}

To construct $F_2$, we move one element from column 1 to column 3, one element from column 1 to column 6, and two elements from column 4 to column 6. We choose the rows to move elements in arbitrarily from the available rows. If we choose rows 7, 1, 2 and 3 respectively, we arrive at the set $F_2$ shown in Figure \ref{plaatjeF2}.

\begin{figure}[h]
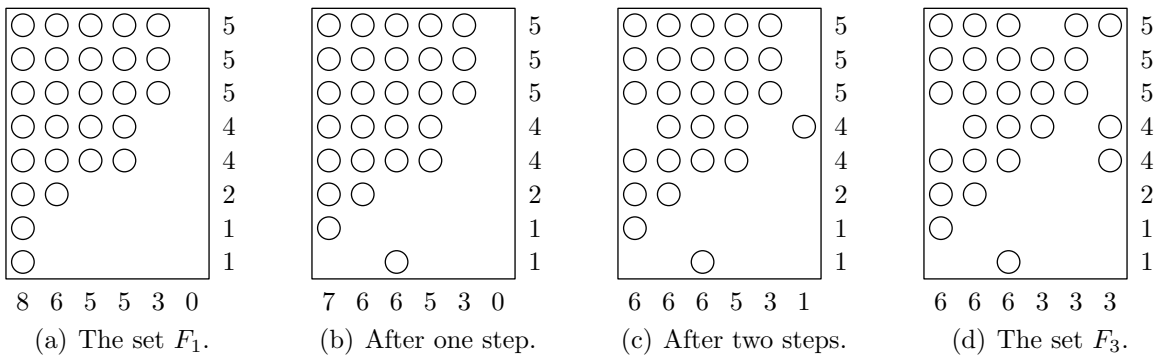

  \begin{center}
    \subfigure[The set $F_1$.]{\label{plaatjeF30}\includegraphics{constructie.1}}
    \qquad
    \subfigure[After one step.]{\label{plaatjeF31}\includegraphics{constructie.5}}
    \qquad
    \subfigure[After two steps.]{\label{plaatjeF32}\includegraphics{constructie.6}}
    \qquad
    \subfigure[The set $F_3$.]{\label{plaatjeF33}\includegraphics{constructie.7}}
  \end{center}
  \caption{The construction of the set $F_3$.}
  \label{plaatjeF3}
\end{figure}

Now we construct the set $F_3$ step-by-step, following the proof of the theorem. We start with $F_1$, shown again in Figure \ref{plaatjeF30}. For the first column pair, we need to move an element from column 1 to column 3. The available rows are 6, 7 and 8. We need only two of them, so let us take $R = \{7, 8\}$. Column 3 is the final column in this column pair, so in this column we need to make sure that we achieve a symmetric difference of at least 2 with $F_2$. We have $(8,3) \not\in F_2$, so we are in case 1 and we pick row 8 for our staircase. Hence we delete the element $(8,1)$ and add the element $(8,3)$. The new situation is shown in Figure \ref{plaatjeF31}.
\\ \\
\\ \\
The next column pair is $(1,6)$. Now column 1 is the final column pair, and all rows except row 8 are available. We are again in case 1 and pick row 4. Figure \ref{plaatjeF32} shows the new situation, after deleting $(4,1)$ and adding $(4,6)$.

Finally, we need to move two elements at once for the column pair $(4,6)$, which occurs twice. Column 6 is the final column, so we need to achieve a symmetric difference of at least 4 with $F_2$ in this column. We also need a symmetric difference of at least 2 in column 4 (case (B)). We have $R = \{1, 2, 3, 5\}$. As $(1,8)$, $(2,8)$ and $(3,8)$ are all elements of $F_2$, we are in case 2. We have $R' = \{1, 2, 3\}$ and we need to find an $l_0 \in R'$ such that $(l_0, 4) \in F_2$. The only possible choice is $l_0 = 1$. We find $R'' = \{1, 5\}$, so we delete $(1,4)$ and $(5,4)$, and we add $(1,6)$ and $(5,6)$. This completes the construction of $F_3$. The resulting set is shown in Figure \ref{plaatjeF33}.

The construction guarantees that the symmetric difference between $F_2$ and $F_3$ is at least $2\alpha +2 = 10$, but we have in fact constructed two sets with symmetric difference 14.

\end{document}